\documentclass[10pt,a4paper]{smfart}

\usepackage[english,francais]{babel}
\usepackage[latin1]{inputenc}
\usepackage{amsmath}
\usepackage{amsfonts}
\usepackage{amssymb}
\usepackage{smfthm}

\newtheorem{propriete}{Propriété}
\newenvironment{hyp}{\begin{propriete}\rm}{\end{propriete}}

\newcommand{\im}{\operatorname{im}}

\newcommand{\emp}[2]{\genfrac{}{}{0pt}{}{#1}{#2}}
\newcommand{\lpro}[1]{\lim\limits_{\emp{\longleftarrow}{#1}}}
\newcommand{\isom}{\stackrel{\sim}{\longrightarrow}}

\newcommand{\spec}{\operatorname{Spec}}
\newcommand{\End}{\operatorname{End}}
\newcommand{\Lie}{\operatorname{Lie}}
\newcommand{\tr}{\operatorname{tr}}
\newcommand{\Hom}{\operatorname{Hom}}
\newcommand{\comp}{\circ}

\newcommand{\chk}[1]{#1^{\vee}}
\newcommand{\point}{{\scriptscriptstyle\bullet}}

\newcommand{\Rep}{\operatorname{Rep}}

\title[]{Comparaison entre cohomologie cristalline et cohomologie étale
$p$-adique sur certaines variétés de Shimura}

\alttitle{Comparison between crystalline cohomology and $p$-adic étale
cohomology on certain Shimura varieties}

\author{Sandra Rozensztajn}
\address{IRMA \\ Université Louis Pasteur \\ 7 rue René-Descartes \\
67084 Strasbourg Cedex \\ France}
\email{rozenszt@math.u-strasbg.fr}

\begin{document}

\begin{abstract}
Soit $X$ un modèle entier en un premier $p$ d'une variété de Shimura de
type PEL, ayant bonne réduction associée à un groupe réductif $G$. On
peut associer aux $\mathbb{Z}_p$-représentations du groupe $G$ deux types
de faisceaux : des cristaux sur la fibre spéciale de $X$, et des systèmes
locaux pour la topologie étale sur la fibre générique. Nous établissons
un théorème de comparaison entre la cohomologie de ces deux types de
faisceaux.
\end{abstract}

\begin{altabstract}
Let $X$ be an integral model at a prime $p$ of a Shimura variety of PEL
type having good reduction, associated
to a reductive group $G$. To $\mathbb{Z}_p$ reprsententations of the group $G$
can be associated two kinds of sheaves : crystals on the special fiber of
$X$, and locally constant étale sheaves on the generic
fiber. We establish a comparison between the cohomology of these two
kinds of sheaves.
\end{altabstract}

\maketitle

\section{Introduction}

Considérons $X$ un modèle entier d'une variété de Shimura de type PEL,
défini sur une extension de $\mathbb{Z}_p$. Cette variété de Shimura correspond
à un groupe réductif $G$, défini sur $\mathbb{Z}_{(p)}$. On peut associer aux
$\mathbb{Z}_{(p)}$-représentations du groupe $G$ différents faisceaux : des
systèmes locaux en $\mathbb{Z}_p$-modules sur la fibre générique de $X$, et des
cristaux sur sa fibre spéciale. Nous établissons une comparaison entre la
cohomologie étale du système local d'une part, et la cohomologie
log-cristalline d'une extension du cristal à une compactification
appropriée de $X$, pour une même représentation $V$ du groupe $G$.
Nous traitons ici le cas des variétés de Shimura unitaires et de celles
de type Siegel. Dans le cas unitaire, nous obtenons un résultat qui tient
compte de la torsion (théorème \ref{principal_torsion}), dans le cas
Siegel les résutats sont moins précis et ne sont valables qu'après
tensorisation par $\mathbb{Q}_p$ (théorème \ref{principal_rationnel}).
L'intérêt de cette comparaison est que nous pouvons obtenir des
renseignements sur le côté cristallin : des techniques de type complexe
BGG, décrites par exemple dans \cite{CF}, chapitre VI,
permettent d'avoir des renseignements sur la filtration de Hodge. On
déduit alors des informations sur le côté étale, vu comme représentation
galoisienne.

\bigskip

La théorie de Hodge $p$-adique nous donne de tels théorèmes de
comparaison entre cohomologie étale $p$-adique et log-cristalline dans le
cas des coefficients constants, pour des schémas propres et possédant
certaines propriétés de lissité. Ces résultats sont dûs entre autres à Tsuji
(\cite{Tsu99}) pour le cas o\`u l'on considère les groupes de cohomologie
après tensorisation par $\mathbb{Q}_p$, et à Tsuji et Breuil pour le cas o\`u
l'on tient compte de la torsion (\cite{Tsu00} et \cite{Bre98}, il y a
alors des restrictions sur le degré des groupes de cohomologie que l'on
peut étudier). Ces théorèmes sont rappelés dans le paragraphe
\ref{comparaison-constant}.

Le principe de notre méthode est de considérer la cohomologie à
coefficients constants de la variété abélienne universelle sur $X$ et de
ses puissances, et d'en déduire la comparaison qui nous intéresse en
découpant les groupes de cohomologie des faisceaux considérés
dans les groupes de cohomologie à coefficients constants à l'aide de
certaines correspondances algébriques.

\bigskip

Le premier problème est que de tels théorèmes de comparaison n'étant
valables que sur des schémas propres, nous devons supposer l'existence
(prouvée dans certains cas seulement) de 
compactifications non seulement de $X$, mais aussi des variétés de
Kuga-Sato, et plus précisément un système projectif de telles
compactifications. La partie \ref{situation} explique précisément dans quelle
situation nous nous pla\,cons, ainsi que les propriétés des
compactifications que nous utilisons.

Le deuxième problème est que toutes les représentations de $G$ ne donnent
pas des faisceaux dont la cohomologie puisse être découpée par des
correspondances algébriques dans la cohomologie de la variété abélienne
universelle.
On détermine dans la partie \ref{representations}
quelles sont les représentations de $G$ qui donnent des faisceaux que
l'on peut atteindre de cette fa\,con, qui sont les seuls pour lesquels nous
obtenons un résultat. Cette partie utilise fortement la structure des
représentations du groupe réductif $G$, ce qui explique que l'on soit
obligé de faire une description cas particulier par cas particulier.

Nous expliquons la construction des faisceaux ainsi que l'action de ces
correspondances algébriques dans la partie \ref{faisceaux}.

Enfin l'énoncé et la preuve du théorème principal occupent la partie
\ref{comparaison}. C'est ici qu'apparaît une différence entre les cas unitaire et
Siegel : en effet le point-clé de la preuve est la compatibilité de
l'action des correspondances algébriques que nous considérons avec les
théorèmes de comparaison à coefficients constants. Le cas Siegel utilise
la compatibilité de cet isomorphisme de comparaison avec les structures
produit sur les groupes de cohomologie, ce qui n'est démontré que dans le
cas rationnel et non dans le cas de torsion.

\section{Les objets considérés}
\label{situation}

\subsection{Variétés de Shimura de type PEL}

\subsubsection{Les données}
\label{donnees-shimura}

On se donne $B$ une $\mathbb{Q}$-algèbre simple finie, munie d'une
involution positive notée $*$ (c'est-à-dire que $\tr_{B/\mathbb{Q}}(xx^*)
> 0$ pour tout $x$ non nul de $B$), $\mathfrak{V}$ un module de type fini
sur $B$, muni d'une forme bilinéaire $(,)$ telle que pour tous $v$ et $w$
dans $\mathfrak{V}$, et tout $b$ dans $B$ on ait $( bv,w ) = ( v,b^*w)$.
On notera $2g$ la dimension de $\mathfrak{V}$ sur $\mathbb{Q}$.

On fixe dans toute la suite un nombre premier $p$, et on fera l'hypothèse
que $p > 2g$. Le rôle de cette hypothèse est expliqué au paragraphe
\ref{unicite-prolongement-connexion}.

On suppose que $B$ est non ramifié en $p$, c'est-à-dire que
$B_{\mathbb{Q}_p}$ est un produit d'algèbres de matrices sur des
extensions non ramifiées de $\mathbb{Q}_p$.

On se donne un $\mathbb{Z}_{(p)}$-ordre $\mathcal{O}_B$ dans $B$ qui
devient un ordre maximal de $B_{\mathbb{Q}_p}$ après tensorisation par
$\mathbb{Z}_p$, et stable par l'involution de $B$.

On se donne aussi $\mathcal{V}$ un $\mathcal{O}_B$-réseau de
$\mathfrak{V}$ autodual.  Le fait que $\mathcal{V}$ soit autodual
implique en particulier que la forme bilinéaire induite sur $\mathcal{V}$
est non dégénérée.

Soit $C$ l'anneau des endomorphismes $B$-linéaires de $\mathfrak{V}$.

On définit le groupe $G$ par $G(R)=\{ g \in (C\otimes R)^*, \exists \mu
\in R^*, \forall v,w \in \mathfrak{V}\otimes R, ( gv,gw ) = \mu ( v,w )
\}$, pour toute $\mathbb{Z}_{(p)}$-algèbre $R$.

On se donne un $\mathbb{R}$-homomorphisme d'algèbres $h : \mathbb{C}
\rightarrow C_{\infty} = C\otimes_{\mathbb{Q}}\mathbb{R}$ tel que
$h(z)^*=h(\bar{z})$, et la forme $( v, h(i)w )$ soit définie positive sur
$\mathfrak{V}_{\infty}=\mathfrak{V}\otimes_{\mathbb{Q}}\mathbb{R}$.  On
associe à $h$ le morphisme $\mu_h : \mathbb{C}^* \to G_{\mathbb{C}}$, qui
définit la filtration de Hodge sur $\mathfrak{V}_{\mathbb{C}}$,
c'est-à-dire la décomposition $\mathfrak{V} = \mathfrak{V}_{z} \oplus
\mathfrak{V}_{1}$, o\`u $\mu_h(z)$ agit par $z$ sur $\mathfrak{V}_{z}$ et
par $1$ sur $\mathfrak{V}_{1}$.

Le corps dual associé à ces données est le corps $E(G,h)$ qui est le
corps de définition de la classe d'isomorphisme de $\mathfrak{V}_{z}$ comme
$B$-représentation. C'est le sous-corps de $\mathbb{C}$ 
engendré par les $\tr(b),b\in B$ agissant sur $\mathfrak{V}_{z}$.

\subsubsection{Deux cas particuliers}

Dans la suite nous nous intéresserons uniquement à deux cas particuliers
: le cas Siegel et le cas unitaire. Le cas Siegel correspond à la
situation o\`u $B$ est réduit à $\mathbb{Q}$. La variété de Shimura
associée est alors la variété modulaire de Siegel. Le cas unitaire
correspond au cas o\`u $B$ est une extension quadratique imaginaire de
$\mathbb{Q}$. La forme alternée $(,)$ est alors la partie imaginaire
d'une forme hermitienne sur $\mathfrak{V}$, qu'on peut voir comme un
$B$-espace vectoriel de dimension moitié.

\subsubsection{Le problème de modules}

On peut associer aux données de Shimura prcédentes un problème de
modules, tel que décrit dans \cite{Kot92}, dont on rappelle ici
l'essentiel.

Fixons ${\mathcal{K}}^p$ un sous-groupe compact ouvert de
$G(\mathbb{A}_f^p)$.  On considère le foncteur des
$\mathcal{O}_{E(G,h)}\otimes\mathbb{Z}_{(p)}$-schémas dans les ensembles,
qui à $S$ associe l'ensemble à équivalence près des quadruplets
$(A,\lambda,i,\eta)$, o\`u $A$ est un schéma abélien $A$ sur $S$, muni
d'une polarisation $\lambda$ première à $p$, et d'une flèche $i :
\mathcal{O}_B \rightarrow \operatorname{End}(A)\otimes\mathbb{Z}_{(p)}$
qui est un morphisme d'algèbres à involution, l'involutin sur
$\End(A)\otimes\mathbb{Z}_{(p)}$ étant l'involution de Rosati donnée par
$\lambda$, et $\eta$ est une structure de niveau. Enfin on suppose que
$\mathcal{O}_B$ agit sur $\Lie(A)$ comme sur $\mathfrak{V}_{z}$,
c'est-à-dire que $\det(b,\Lie(A)) = \det(b,\mathfrak{V}_{z})$ pour tout
$b \in \mathcal{O}_B$.  La structure de niveau consiste en ce qui suit :
on considère le $\mathbb{A}_f^p$-module de Tate de $A$, c'est un
$\mathbb{A}_f^p$-faisceau lisse sur $S$. Soit $s$ un point géométrique de
$S$, une structure de niveau consiste en une ${\mathcal{K}}^p$-orbite
$\eta$ d'isomorphismes $\mathfrak{V}_{\mathbb{A}_f^p} \rightarrow
H_1(A_s,\mathbb{A}_f^p)$ de $B$-modules munis d'une forme alternée, et
qui soit fixée par $\pi_1(S,s)$.

Deux quadruplets $(A,\lambda,i,\eta)$ et $(A',\lambda',i',\eta')$ sont
dit équivalents s'il existe une isogénie première à $p$ de $A$ vers $A'$,
commutant à l'action de $\mathcal{O}_B$, transformant $\eta$ en $\eta'$, et
$\lambda$ en un multiple scalaire (dans $\mathbb{Z}_{(p)}^*$) de $\lambda'$.

Ce foncteur est représentable, par un schéma quasi-projectif et lisse
$\mathcal{M}$ sur $\mathcal{O}_E\otimes\mathbb{Z}_{(p)}$ pourvu
que l'on choisisse ${\mathcal{K}}^p$ suffisamment petit.

Notons $\mathcal{A}$ le schéma abélien universel sur
$\mathcal{M}$.

Pour les cas Siegel et unitaire, 
on a le résultat suivant (lemme 7.2 de \cite{Kot92}) :
\begin{lemm}
\label{torseur-anneau}
$\mathcal{O}_{\mathcal{M}}\otimes \mathcal{V}$ et
$\mathcal{H}^1(\mathcal{A}/\mathcal{M}\chk{)}$ sont
localement isomorphes comme $\mathcal{O}_B$-modules munis d'une forme alternée.
\end{lemm}

\subsubsection{La situation géométrique considérée}

On obtient alors la situation suivante : Notons $K$ le complété en une
place $v |p$ de $E(G,h)$, et $\mathcal{O}$ son anneau des entiers. Notons
$\mathcal{O}_n=\mathcal{O}/\varpi^{n+1}$, o\`u $\varpi$ est une
uniformisante de $\mathcal{O}$. C'est l'anneau des vecteurs de Witt de
longueur $n$ sur le corps résiduel de $\mathcal{O}$ puisque $p$ est non
ramifié dans $E(G,h)$.  Posons $S = \spec \mathcal{O}$, et $S_n = \spec
\mathcal{O}_n$.  D'une fa\,con générale, on notera avec un indice $n$ la
réduction d'un $\mathcal{O}$-schéma modulo $\varpi^{n+1}$

On notera $X = \mathcal{M}_{\mathcal{O}}$, c'est donc un schéma lisse sur
$S$, muni d'un schéma abélien $\mathcal{A}$, provenant du schéma abélien
universel sur la variété de Shimura. De plus, on a un morphisme
$\mathcal{O}_B \to \End(A)\otimes_{\mathbb{Z}}\mathbb{Z}_{(p)}$.  On note
$f : \mathcal{A} \to X$, et $f_s : \mathcal{A}^s \to X$ les morphismes
structuraux.

\subsection{Existence de compactifications}
\label{hypotheses}

Nous aurons besoin d'utiliser aussi des compatifications de $X$, ainsi
que du schéma abélien universel $\mathcal{A}$ et de ses puissances.  Ces
compactifications ont été décrites en détail dans le cas Siegel
(\cite{CF}), pour le cas unitaire la construction détaillée n'est écrite
que pour $GU(2,1)$, c'est-à-dire les surfaces modulaires de Picard (voir
\cite{Lar92} pour la compactification de la base et \cite{Roz06} pour
celle du schéma abélien universel), même si leur existence dans le cas
unitaire général ne pose pas de problèmes.  Nous résumons dans ce
paragraphe les seules propriétés de ces compactifications que nous
utilisons.

\subsubsection{Compactifications de la base}

Nous avons besoin tout d'abord de compactification du modèle entier de la
variété de Shimura. En considérant des compactifications toroïdales
(décrites dans \cite{CF}, chapitre IV pour le cas Siegel, et dans
\cite{Lar92} pour le cas de $GU(2,1)$), on obtient l'existence d'un
schéma $\overline{X}$ vérifiant la propriété suivante :

\begin{hyp}
\label{hyp1}
il existe un schéma $\overline{X}$ propre et lisse sur $S$, contenant $X$
comme ouvert dense, tel que le complémentaire de $X$ dans $\overline{X}$
est un diviseur à croisements normaux relatifs.
\end{hyp}

On fixera dans la suite une fois pour toute une telle compactification
$\overline{X}$. On peut remarquer que les résultats ne dépendent en fait pas
du choix de $\overline{X}$ parmi l'ensemble des compactifications toroïdales :
deux compactifications toroïdales sont toujours comparables, au sens
o\`u il existe une troisième qui les domine toutes les deux, ce qui permet de
voir que les groupes de cohomologie décrits en \ref{coho-cristaux} ne 
dépendent pas de ce choix de compactification.

\subsubsection{Compactifications du schéma abélien universel}

On se donne $\overline{X}$ propre lisse sur $S$, muni d'un diviseur à
croisements normaux relatifs $D$, et on note $X$ l'ouvert complémentaire.
Pour chaque $s$, on appelle bonne compactification de $\mathcal{A}^s$ une
compactification $\overline{\mathcal{A}^s}$ de $\mathcal{A}^s$, telle que
$f_s^{-1}(\overline{X} \setminus X)$ est un diviseur à croisements
normaux relatifs. 

Considérons la des compactifications toroïdales « lisses » des
$\mathcal{A}^s$ (voir \cite{CF} dans le cas Siegel, \cite{Roz06} dans le
cas de $GU(2,1)$), on obtient pour tout $s \geq 1$, une famille de bonnes
compactifications $\overline{\mathcal{A}^s}$ de $\mathcal{A}^s$, vérifiant les deux propriétés
suivantes :
\begin{hyp}
\label{hyp2}
\begin{enumerate}
\item Pour toute isogénie $u$ de $\mathcal{A}^s$, il existe deux bonnes
compactifications $\overline{\mathcal{A}^s}_1$ et $\overline{\mathcal{A}^s}_2$, et un morphisme
$\overline{\mathcal{A}^s}_1 \to \overline{\mathcal{A}^s}_2$ prolongeant $u$. 
\item \'Etant donné deux bonnes
compactifications $\overline{\mathcal{A}^s}_1$ et $\overline{\mathcal{A}^s}_2$, il en existe une
troisième $\overline{\mathcal{A}^s}_3$ et des $\overline{X}$-morphismes 
$\overline{\mathcal{A}^s}_3 \to  \overline{\mathcal{A}^s}_1$ et $\overline{\mathcal{A}^s}_3 \to  \overline{\mathcal{A}^s}_2$
induisant l'identité sur $\mathcal{A}^s$.
\end{enumerate}
\end{hyp}

Dans le cas Siegel, nous utiliserons encore une propriété supplémentaire
de la famille des compactifications toroïdales :
\begin{hyp}
\label{hyp3}
Si $\mathcal{L}$ est un faisceau symétrique sur $\mathcal{A}^s$, il existe
des entiers $a$ et $b$, et une compactification de la famille
$\overline{\mathcal{A}^s}$ tels que le faisceau $(\mathcal{O}(2)^{\otimes
a}\otimes \mathcal{L})^{\otimes b}$ se prolonge en un faisceau sur
$\overline{\mathcal{A}^s}$. De plus, on peut choisir $a$ et $b$ premiers à
$p$.
\end{hyp}

Cette construction fait l'objet du chapitre VI du livre \cite{CF}. Elle
n'est détaillée que pour le cas o\`u le faisceau symétrique ample considéré
est $\mathcal{O}(2)$, mais cela s'adapte au cas d'un faisceau symétrique ample
quelconque, pour lequel on prendra donc \mbox{$\mathcal{O}(2)^{\otimes a}\otimes
\mathcal{L}$}, avec $a$ assez grand pour que le faisceau soit ample.

\section{Représentations de $G$}
\label{representations}

\subsection{$\mathbb{Z}_{(p)}$-représentations}

On note $\Rep(G)$ la catégorie des représentations de $G$ sur un
$\mathbb{Z}_{(p)}$-module libre de type fini, et $\Rep_{\mathbb{Q}}(G)$
celles des représentations de $G$ sur un $\mathbb{Q}$-espace vectoriel de
dimension finie.

Notons $\mathcal{V}_0 \in \Rep(G)$ la duale de la représentation standard
de $G$, c'est-à-dire la duale du réseau $\mathcal{V}$ défini au
paragraphe \ref{donnees-shimura}.

\subsection{L'algèbre des endomorphismes}
\label{algebre_endo}

Les représentations de $G$ de la forme 
$\wedge^{\point}\mathcal{V}_0^s$, pour $s \geq 1$, jouent un rôle
particulier : les faisceaux que nous allons leur associer dans la section
\ref{faisceaux} ont une interprétation géométrique.  Nous définissons
une sous-algèbre de l'algèbre $\End(\wedge^{\point}\mathcal{V}_0^s)$ 
des endomorphismes $G$-linéaires de
$\wedge^{\point}\mathcal{V}_0^s$, formé de morphismes ayant aussi une
interprétation géométrique qui sera décrite dans le paragraphe
\ref{traduction_geometrique}.  L'objectif est de pouvoir découper
dans les $\wedge^{\point}\mathcal{V}_0^s$ des représentations
irréductibles de $G$ à l'aide de cette algèbre d'endomorphismes, en
s'inspirant des constructions de Weyl.

\smallskip

Dans le cas unitaire, on définit pour tout $s \geq 1$ une
sous-$\mathbb{Z}_{(p)}$-algèbre $E(A)_s$ de
$\End(\wedge^{\point}\mathcal{V}_0^s)$.  C'est l'algèbre engendrée par
l'action de $M_s(\mathbb{Z})$ sur $\mathcal{V}_0^s$, muni de la
multiplication opposée, et par l'action de $\mathcal{O}_B$ sur
$\mathcal{V}_0$.

\smallskip

Dans le cas Siegel, on note $E(C)_s$ la sous-$\mathbb{Z}_{(p)}$-algèbre
de $\End(\wedge^{\point}\mathcal{V}_0^s)$ engendrée par l'action de
$M_s(\mathbb{Z})$ sur $\mathcal{V}_0^s$, et par les opérations suivantes.
On note $\mathbb{Z}_{(p)}(1)$ la représentation du groupe symplectique
correspondant à l'action du groupe sur $\mathbb{Z}_{(p)}$ par le
multiplicateur. Observons que $\mathcal{V}_0 = \mathcal{V}(-1)$.

On note $u \in \wedge^2\mathcal{V}_0^2(1)$ l'élément provenant de la
forme bilinéaire sur $\mathfrak{V}$, et pour tous $1\leq i < j \leq s$, on note
$u_{i,j}$ l'image de $u$ par l'application $\wedge^2\mathcal{V}_0^2(1)
\to \wedge^2\mathcal{V}_0^s(1)$ induite par l'application
$\mathcal{V}_0^2 \to \mathcal{V}_0^s$ consistant à placer les deux
facteurs $\mathcal{V}_0$ aux places $i$ et $j$.

Le cup-produit par $u_{i,j}$ définit une application
$\varphi_{i,j} : \wedge^{\point}\mathcal{V}_0^s(-1) \to
\wedge^{\point}\mathcal{V}_0^s$ qui envoie chaque 
$\wedge^{k}\mathcal{V}_0^s(-1)$ dans
$\wedge^{k+2}\mathcal{V}_0^s$.

L'application duale de $\varphi_{i,j}$ permet de définir une application 
$\psi_{i,j} : \wedge^{\point}\mathcal{V}_0^s \to
\wedge^{\point}\mathcal{V}_0^s(-1)$.

Enfin on note $\theta_{i,j} = \varphi_{i,j}\comp\psi_{i,j}$, qui est donc un
endomorphisme de $\wedge^{\point}\mathcal{V}_0^s$.

On définit alors $E(C)_s$ comme l'algèbre engendrée par
l'action de $M_s(\mathbb{Z})$ et les $\theta_{i,j}$, $1\leq i < j \leq s$.

\smallskip

On notera $E_s$ pour désigner indifféremment $E(A)_s$ et $E(C)_s$.
On dira que $u \in E_s$ est un projecteur homogène (de degré $t$) si son
image est contenue dans $\wedge^t\mathcal{V}_0^s \subset \wedge^{\point}\mathcal{V}_0^s$.

\smallskip

On donne des définitions similaires pour les éléments de $E_s\otimes\mathbb{Q}$
agissant sur \mbox{$\wedge^{\point}(\mathcal{V}_0\otimes\mathbb{Q})^s$}.

\subsection{Représentations atteignables}
\label{atteignable}

On note $\Rep^a(G)$ la sous-catégorie de $\Rep(G)$ formée des
représentations isomorphes à une somme directe de représentations de $G$
de la forme $\im q$, o\`u $q$ est un projecteur homogène de $E_s$
agissant sur un $\wedge^{\point}\mathcal{V}_0^s$. Si $V\in \Rep^a(G)$, on note
$t(V)$ le plus grand degré des projecteurs homogènes qui apparaissent
dans la définition de $V$.

On définit de fa\,con similaire $\Rep^a_{\mathbb{Q}}(G)$.

Comme nous n'obtenons des résultats que pour les représentations de $G$
qui sont dans $\Rep^a(G)$, il va s'agir de voir que cette sous-catégorie
n'est pas trop petite, et qu'il n'est donc pas trop restrictif de s'y
limiter. C'est l'objet des paragraphes suivants.

\subsection{Poids $p$-petits}

Supposons notre groupe réductif $G$ déployé sur un certain corps $E$. Les
représentations irréductibles de $G$ sur $E$ sont paramétrées par
l'ensemble des poids dominants, une fois fixé un tore maximal et un
système de racines positives. 
Si $a$ est un poids dominant, on note $V_E(a)$ la représentation
irréductible de $G$ sur $E$ de plus haut poids $a$.

On définit comme dans \cite{Jan87}, II.3.15 ce qu'est un poids dominant
$p$-petit.  La propriété qui nous intéresse ici est la propriété suivante
(voir \cite{PolTil02}, 1.9) : si le poids $a$ est $p$-petit, alors il
existe à homothétie près un unique réseau dans $V_E(a)$ qui est stable
sous l'action de $\text{Dist}(G)$, l'algèbre des distributions de $G$ sur
$\mathcal{O}_{E,(v)}$, pour $v$ une place de $E$ divisant $p$.  Cela a
donc un sens de définir $V(a)$ comme la représentation irréductible de
$G$ sur $\mathcal{O}_{E,(v)}$ de plus haut poids $a$.

\subsection{Description de $\Rep^a(G)$ dans le cas unitaire}

On se place ici dans le cas unitaire.  Le groupe $G$ est alors un groupe
unitaire relatif à un corps $E$ quadratique imaginaire, qui correspond à
l'algèbre $B$ du paragraphe \ref{donnees-shimura}, donc $G$ est de la forme
$GU(g)$, et il est déployé sur $E$. On a donc $G_E \isom GL(g)_E \otimes
\mathbb{G}_{m,E}$.

L'ensemble des représentations sur $E$ irréductibles de $G$ est paramétré
par les $g+1$-uplets $(a_1,\dots,a_g;c)$ d'entiers, avec $a_1 \geq \dots
\geq a_g$ et $\sum a_i = c \pmod 2$, une fois choisi un isomorphisme
entre $G_E$ et $GL(g)_E \times \mathbb{G}_{m,E}$. Notons $i$ et $j$ les
deux morphismes de $E$ dans $E$, le choix de l'isomorphisme revient à en
privilégier un des deux.

La représentation $\mathcal{V}_0\otimes E$ correspond à la somme de deux
représentations irréductibles $V_1$ et $V_2$, $V_1$ de plus haut poids
$(1,0,\dots,0;1)$ et $V_2$ de plus haut poids $(0,\dots,0,-1;1)$. $V_1$
est l'espace propre associé à la valeur propre $i(x)$ de l'endomorphisme
$u(x)$, pour tout $x$ dans $E$, et $V_2$ est l'espace propre associé à la
valeur propre $j(x)$.

\subsubsection{Description de $\Rep_{\mathbb{Q}}^a(G)$}
\label{rep-Q}

Notons $\Rep_E^a(G)$ l'ensemble des $V\otimes_{\mathbb{Q}}E$, o\`u $V \in
\Rep_{\mathbb{Q}}^a(G)$, et $V_0 = \mathcal{V}_0\otimes E$.

\begin{prop}
$\text{Rep}^a_{E}(G)$ contient toutes les représentations qui sont de la
forme $V(a)\oplus V(a^*)$, o\`u $V(a)$ est la représentation irréductible
de plus haut poids $(a_1,\dots,a_g;c)$ avec $a_g \geq 0$ et $c = \sum
a_i = s$, et $a^*$ est le  poids $(-a_g,\dots,-a_1;c)$.
\end{prop}

\begin{lemm}
\label{schur}
Soit $a = (a_1,\dots,a_g;c)$ un poids dominant de $G$, tel que $a_g \geq
0$
et $c = \sum a_i$, et $V(a)$ la représentation irréductible associée.
Alors il existe un élément $C_a$ de $\mathbb{Q}\mathfrak{S}_s$ (o\`u $s = \sum
a_i$) tel que $V(a) = C_aV_1^{\otimes s}$ et $V(a^*) = C_aV_2^{\otimes
s}$.
\end{lemm}

\begin{proof}
Regardons d'abord $V(a)$ comme une représentation de $GL_g$, en oubliant
l'action du multiplicateur. Comme expliqué dans \cite{FulHar}, 15.5, il
existe un $C_a \in \mathbb{Q}\mathfrak{S}_s$ idempotent tel que
$V(a) = C_aV_1^{\otimes s}$, $V(a)$ et $V_1$ étant vues toutes deux comme
des
représentations de $GL_g$. Il faut voir ensuite que l'égalité tient aussi
comme représentations de $GL(g)_E\times \mathbb{G}_{m,E}$, donc que le
multiplicateur agit de la même fa\,con sur les deux. Or il agit par $x
\mapsto
x^c$ sur $V(a)$, et par $x \mapsto x^s$ sur $V_1^{\otimes s}$, et on a $s =
c$.
\end{proof}

\begin{lemm}
\label{tens}
Soit $s \geq 0$. Il existe un projecteur $q$ dans $E(A)_s\otimes\mathbb{Q}$ commutant à
l'action du groupe des permutations $\mathfrak{S}_s$ tel
que l'image de $q$ agissant sur $\wedge^{\point}V_0^s$
est $V_0^{\otimes s}$.
\end{lemm}

\begin{proof}
$\wedge^{\point}V_0^s = \oplus_{0 \leq i_1 \leq 2g, \dots, 0 \leq i_s
\leq 2g} \wedge^{i_1}V_0\otimes\dots\otimes\wedge^{i_s}V_0$.

Considérons un entier $m$ non nul, $v_j$ la matrice diagonale
$[(1,\dots,1,m,1,\dots)]$ avec un $m$ en $j$-ème position. L'espace
propre correspondant à la valeur propre $m$ est
la somme des termes pour lesquels $i_j = 1$. Soit $p_j$ le
projecteur sur cet espace propre. Les $p_j$ commutent, leur produit est
donc un
projecteur $q$ sur l'intersection des images, c'est-à-dire les termes
pour lesquels chaque $i_j$ est égal à $1$, c'est-à-dire $V_0^{\otimes
s}$. De plus, $q$ commute bien à l'action de $\mathfrak{S}_s$.
\end{proof}

Enfin on utilise que $V_0 = V_1 \oplus V_2$, $V_0^{\otimes s}$ est donc
égal à une somme de termes de la forme $V_1^{\otimes x}\otimes
V_2^{\otimes y}$ avec $x+y = s$.

\begin{lemm}
\label{E}
Il existe un élément $q'$ dans le centre de $E(A)_s\otimes\mathbb{Q}$ dont la restriction de
l'action à $V_0^{\otimes s}$ est un projecteur sur $V_1^{\otimes s}\oplus
V_2^{\otimes s}$.
\end{lemm}

\begin{proof}
Fixons $z\in \mathcal{O}_E$. $z$ agit par $i(z)$ sur $V_1$ et par $j(z)$ sur
$V_2$.
Notons $i(z) = a$ et $j(z) = b$. Choisissons $z$ de sorte que les
$a^xb^y$ soient tous distincts. Alors
$V_1^{\otimes x}\otimes V_2^{\otimes y}$ est (dans $V_0^{\otimes s}$)
l'espace propre associé à la valeur propre $a^xb^y$.

Soit $P$ le polynôme $\Pi_{r+t=s}(X-a^rb^t)$. Il est à coefficients
entiers, et c'est le polynôme minimal de l'action de $u(z)$ sur
$V_0^{\otimes
s}$. Soit $Q = \Pi_{r+t=s, r\neq 0, t\neq 0}(X-a^rb^t)$.
Alors $P = QT$ o\`u $T = (X-a^s)(X-b^s)$, et $Q$ et $T$ sont premiers
entre eux. Il existe donc des polynômes $U$ et $V$ (à coefficients
rationnels), tels que $UQ + VT = 1$. Alors l'action de $u(1-UQ)(z)$ sur
$V_0^{\otimes s}$ est un projecteur sur $V_1^{\otimes s}\oplus
V_2^{\otimes s}$, qu'on note $q'$.
\end{proof}

On peut maintenant prouver la proposition : soit $a$ comme dans l'énoncé,
et $s = c = \sum a_i$. Posons $P = C_aq'q$, alors $P\wedge^{\point}V_0^s$
est la représentation $V(a)\oplus V(a^*)$. En effet, notons que $C_a$,
$q'$ et $q$ commutent par construction, donc $P$ est un projecteur. On a
$q\wedge^{\point}V_0^s = V_0^{\otimes s}$, $q'q\wedge^{\point}V_0^s
= V_1^{\otimes s}\oplus V_2^{\otimes s}$,
$C_aq'q\wedge^{\point}V_0^s = C_aV_1^{\otimes s}\oplus C_aV_2^{\otimes
s}$.
Or $C_aV_1^{\otimes s} = V(a)$, et $C_aV_2^{\otimes s} = V(a^*)$.

\bigskip

Une fois décrites les représentations qui sont dans $\Rep^a_E(G)$, il
faut maintenant retrouver quelle est la $\mathbb{Q}$-forme de ces
représentations qui est dans $\Rep^a_{\mathbb{Q}}(G)$. $V_1$ et $V_2$
sont naturellement isomorphes, et $\text{Gal}(E/\mathbb{Q})$ agit sur
$V_0 = V_1 \oplus V_2$ par $(x,y) \mapsto (\bar{y}, \bar{x})$. Son action
sur $V(a) \oplus V(a^*)$ peut être décrite par la même formule, ce qui
nous permet d'obtenir les représentations qui sont dans
$\Rep^a_{\mathbb{Q}}(G)$.

\subsubsection{Description de $\Rep^a(G)$}

\begin{prop}
$\Rep^a(G)$ contient toutes les représentations qui sont de la
forme $(V(a)\oplus V(a^*))^{\text{Gal}(E/\mathbb{Q})}$, o\`u $V(a)$ est la représentation irréductible
de plus haut poids $(a_1,\dots,a_g;c)$ avec $a_g \geq 0$ et $c = \sum
a_i = s$, et $a^*$ est le  poids $(-a_g,\dots,-a_1;c)$, $a$ et $a^*$ sont
$p$-petits, et $2g < p$, et $\text{Gal}(E/\mathbb{Q})$ agit sur $(V(a)\oplus
V(a^*))$ comme décrit au paragraphe précédent.
\end{prop}

Il suffit de voir que si les conditions données sont vérifiées, on peut
prendre des dénominateurs premiers à $p$ dans les lemmes du paragraphe
\ref{rep-Q}.

\begin{lemm}
On se place comme dans le lemme \ref{schur}. Alors si $\sum a_i < p$,
$C_a$ est dans $\mathbb{Z}_{(p)}\mathfrak{S}_s$.
\end{lemm}

\begin{proof}
$C_a$ est de la forme $(1/n)C'_a$, o\`u $C'_a$ est dans $\mathbb{Z}\mathfrak{S}_s$,
et
$n$ est l'entier tel que ${C'_a}^2 = nC'_a$. Or $n$
divise $s!$ (voir \cite{FulHar}, 4.2), donc $1/n \in \mathbb{Z}_{(p)}$ si $s<p$.
\end{proof}

\begin{lemm}
Dans le cadre du lemme \ref{tens}, on peut choisir $q$ dans $E(A)_s$ dès que
$p > 2g$.
\end{lemm}

\begin{proof}
Lorsque on écrit $p_j$ comme un polynôme en $v_j$, les déno\-mi\-na\-teurs qui
apparaissent sont les différences entre les valeurs propres de $v_j$, qui
sont les $m^i$ pour $0 \leq i \leq 2g$. Si $2g < p$, on peut prendre un
$m$ dont l'image dans $\mathbb{Z}/p\mathbb{Z}$ est un générateur de $\mathbb{Z}/p\mathbb{Z}^*$, de sorte
que les $m^i - m^{i'}$ sont tous premiers à $p$.
\end{proof}

\begin{lemm}
Dans le lemme \ref{E}, on peut prendre $q'$ dans $E(A)_s$ dès que $p > s$.
\end{lemm}

\begin{proof}
\'Ecrivons donc $UQ + VT = c$, avec $U$ et $V$ à coefficients entiers et
$c$ entiers, et étudions les facteurs premiers de $c$.  On obtient $c =
Q(a^s)Q(b^s)$.  Il s'agit donc de trouver $z\in \mathcal{O}_E$ tel que $c$ soit
premier à $p$ (et que aucun des $a^rb^t, r>0,t>0, r+t=s$ ne soit égal à
$a^s$ ou à $b^s$).

On a $Q(a^s)=a^{s(s-1)/2}\Pi_{1\leq t \leq s-1}(a^t - b^t)$, et
$Q(b^s)=b^{s(s-1)/2}\Pi_{1\leq t \leq s-1}(b^t-a^t)$.  Supposons d'abord
que $p$ est inerte dans $E$. Alors $\mathcal{O}_E/p$ est isomorphe à
$\mathbb{F}_{p^2}$. La conjugaison dans $\mathcal{O}_E$ se traduit par $x\mapsto
x^p$ dans $\mathcal{O}_E/p$.  Choisissons donc $x$ un générateur de
$\mathbb{F}_{p^2}^*$, alors un $z$ relevant $x$ convient.

Supposons maintenant $p$ décomposé dans $E$.  Alors $\mathcal{O}_E/p$ est égal à
$\mathbb{F}_p \times \mathbb{F}_p$, et la conjugaison dans $\mathcal{O}_E$ échange
les deux facteurs dans $\mathcal{O}_E/p$. Choisissons un $x$ dans $\mathbb{F}_p^*$
tel que $x^i \neq 1$ pour tout $i$ entre $1$ et $s$, et prenons $u$ et
$v$ dans $\mathbb{F}_p^*$ tels que $u/v=x$. Alors si $z$ est un
relèvement de $(u,v)$, il convient.
\end{proof}

\subsection{Description de $\Rep^a(G)$ dans le cas Siegel}
Dans le cas Siegel, la description de $\Rep^a(G)$ est faite dans
l'article \cite{MokTil02}, 5.1. On obtient toutes les représentations de
plus haut poids $p$-petit, à l'action du centre près.

\section{Les faisceaux}
\label{faisceaux}

\subsection{Constructions fonctorielles}

\subsubsection{Cas étale}
\label{fibres-etale}

Soit $x$ un point géométrique de $X_K$.  \`A chaque représentation du
groupe fondamental de la variété $\pi_1(X_K,x)$ correspond un système
local sur $X_K$.  Considérons le faisceau constant $\mathbb{Z}_p$ sur
$\mathcal{A}$, et $F = R^1f_{K*}\mathbb{Z}_p(1)$, o\`u $f_K :
\mathcal{A}_K \to X_K$ est le morphisme structural. $F$ correspond à la
représentation standard de $G$ sur $\mathbb{Z}_p$, autrement dit à un
morphisme $\pi_1(X_K,x) \to G(\mathbb{Z}_p)$. On peut donc associer par
composition un système local à toute représentation sur $\mathbb{Z}_p$ de
$G$, et ceci de fa\,con fonctorielle. On note $\mathbb{F}(V)$ le système
local associé à la représentation $V$. On définit de même le foncteur
$\mathbb{F}_n(V)$, qui à $V$ associe un fibré en
$\mathbb{Z}/p^n\mathbb{Z}$-modules, vérifiant
$\mathbb{F}(V)\otimes_{\mathbb{Z}_p}\mathbb{Z}/p^n\mathbb{Z} =
\mathbb{F}_n(V)$.

On observe que $\mathbb{F}(\mathcal{V}_0) = R^1f_*\mathbb{Z}_p$, et plus
généralement $\mathbb{F}(\wedge^t\mathcal{V}_0^s) = 
R^tf_{s,K*}\mathbb{Z}_p$, o\`u $f_{s,K}$ est le morphisme structural $\mathcal{A}^s_K \to X_K$.

\subsubsection{Cas des fibrés à connexion}

\begin{prop}
Il existe un foncteur $\mathcal{F}$ de $\text{Rep}(G)$ dans
l'ensemble des $\mathcal{O}_X$-modules à connexion intégrable sur $X$, et pour
tout $n$ un foncteur $\mathcal{F}_n$ de $\text{Rep}(G)$ dans
l'ensemble des $\mathcal{O}_{X_n}$-modules à connexion intégrable sur $X_n$, ces
deux foncteurs étant compatibles. 
\end{prop}

Ici compatible, signifie que $\mathcal{F}_n(V/\varpi^{n+1}) =
\mathcal{F}(V)/\varpi^{n+1}$.  On ne va faire la construction que sur
${\mathcal{O}}_X$, la construction sur ${\mathcal{O}}_{X_n}$ s'obtenant par des méthodes
similaires.

Notons  $\mathcal{H}_1(\mathcal{A})=
\chk{(R^1f_*\Omega^{\point}_{\mathcal{A}/X})}$.  Introduisons
$\mathcal{T}=\mathit{Isom}(\mathcal{O}_X\otimes
\mathcal{V},\mathcal{H}_1(\mathcal{A}))$, les isomorphismes devant
respecter la structure de $B$-module et la forme alternée à une constante
près. C'est un torseur sur $X$ sous l'action (à droite) de $G$, en effet
les deux faisceaux en question sont localement isomorphes, comme expliqué
dans le lemme \ref{torseur-anneau}.

Soit $V \in \text{Rep}(G)$, on note $\mathcal{F}(V)$ le faisceau des
sections du fibré $\mathcal{T}\times^GV$. C'est un faisceau de
$\mathcal{O}_X$-modules quasi-cohérent. De même un morphisme entre
représentations se transforme en morphisme entre faisceaux.

Ce fibré est muni d'une connexion qui provient de la connexion de
Gauss-Manin sur $\mathcal{H}_1(\mathcal{A})$.

Notons $\mathcal{H}^i(\mathcal{A}^s) =
R^if_{s*}\Omega^{\point}_{\mathcal{A}^s/X}$, on a alors
$\mathcal{F}(\wedge^t\mathcal{V}_0^s) = \mathcal{H}^t(\mathcal{A}^s)$.
De même on notera $\mathcal{H}^i(\mathcal{A}^s_n) = R^if_{s*}\Omega^{\point}_{\mathcal{A}^s_n/X_n}$.

\subsection{Action de $E_s$}

\subsubsection{Traduction géométrique}
\label{traduction_geometrique}

Comme $E_s$ est une sous-algèbre de
$\End(\wedge^{\point}\mathcal{V}_0^s)$, par fonctorialité de $\mathbb{F}$
et $\mathcal{F}$, on a donc aussi des morphismes de
$\mathbb{Z}_{(p)}$-algèbres $E_s \stackrel{a_{\text{ét}}}{\to} \End
(R^{\point}f_{s,*}\mathbb{Z}_p)$ et $E_s \stackrel{a_{\text{cris}}}{\to}
\End (\mathcal{H}^{\point}(\mathcal{A}^s))$.  On va donner une
interprétation géométrique de ces deux morphismes.

Soit $\mathcal{G}\subset E_s$ la partie formée des éléments suivants :
les matrices de déterminant non nul, dans le cas unitaire les éléments
non nuls de $\mathcal{O}_B$, dans le cas Siegel les opérations $\theta_{i,j},
1\leq i <j \leq s$ définies au paragraphe \ref{algebre_endo}.
L'ensemble $\mathcal{G}$, qu'on appellera ensemble des
éléments géométriques de $E_s$, engendre $E_s$ comme $\mathbb{Z}_{(p)}$-algèbre.  

Soit $u\in \mathcal{G}$ qui provient d'une matrice ou, dans le cas
unitaire, d'un élément de $\mathcal{O}_B$. Alors $u$ provient d'un
élément de $\End(\mathcal{V}_0^s)$, qui agit naturellement sur
$\mathcal{A}^s/X$, donc sur $R^{\point}f_{s,*}\mathbb{Z}_p$ et
$\mathcal{H}^{\point}(\mathcal{A}^s)$ par $a_{\text{ét}}(u)$ et
$a_{\text{cris}}(u)$ respectivement.

Soit $\mathcal{P}$ le faisceau de Poincaré sur
$\mathcal{A}\times\mathcal{A}$, pour $1\leq i<j\leq s$ on note
$\mathcal{P}_{i,j}$ le faisceau sur $\mathcal{A}^s$ obtenu en tirant
$\mathcal{P}$ par le morphisme de projection sur les $i$-ièmes et
$j$-ièmes facteurs $\mathcal{A}^s \to \mathcal{A}\times\mathcal{A}$.
Alors l'opération consistant à faire le cup-produit par la première
classe de Chern de $\mathcal{P}_{i,j}$ correspond à
$a_{\text{ét}}(\psi_{i,j})$ et $a_{\text{cris}}(\psi_{i,j})$, l'opération
duale correspond à $a_{\text{ét}}(\varphi_{i,j})$ et
$a_{\text{cris}}(\varphi_{i,j})$, comme expliqué dans \cite{MokTil02}.

\smallskip

Par fonctorialité des constructions précédentes, on a, pour une
représentation $V$ de la forme $V = q(\wedge^{\point}\mathcal{V}_0^s)$,
$q$ étant un projecteur de $E_s$ : $\mathbb{F}(V) =
a_{\text{ét}}(q)R^{\point}f_{s,K*}\mathbb{Z}_p$, et $\mathcal{F}(V) =
a_{\text{cris}}(q)\mathcal{H}^{\point}(\mathcal{A}^s)$.

On notera encore $a_{\text{cris}}$ et $a_{\text{ét}}$ les morphismes
naturels de $E_s$ vers $\End(\mathcal{H}^{\point}(\mathcal{A}^s_n))$ et
$\End(R^{\point}f_{s*}\mathbb{Z}/p^n\mathbb{Z})$ respectivement.

\subsubsection{Conséquence sur les faisceaux à connexion}

\begin{lemm}
Pour tout $V \in \text{Rep}^{a}(G)$
la connexion sur $\mathcal{F}(V)$  et sur $\mathcal{F}_n(V)$
est quasi-nilpotente.
\end{lemm}

\begin{proof}
Notons que les opérations élémentaires commutent à la connexion sur
$\mathcal{H}^{\point}(\mathcal{A}^s)$ induite par la connexion de
Gauss-Manin, de sorte que, en reprenant les notations précédentes,
$\mathcal{F}(V)$ est stable par la connexion de
$\mathcal{H}^{\point}(\mathcal{A}^s)$. Comme la connexion de Gauss-Manin
sur $\mathcal{H}^{\point}(\mathcal{A}^s)$ est quasi-nilpotente, c'est
aussi le cas pour la connexion sur $\mathcal{F}(V)$.  
\end{proof}

Comme $X$ est lisse sur $S$, chaque $X_n$ est un relèvement de $X_0$ qui
est lisse sur $S_n$. Les faisceaux $\mathcal{F}_n(V)$, qui sont des
$\mathcal{O}_{X_n}$-modules cohérents munis d'une connexion intégrable et
quasi-nilpotente, définissent donc des cristaux sur
$(X_0/S_n)_{\text{cris}}$, ainsi que dans $(X_m/S_n)_{\text{cris}}$ pour
tout $m \leq n$. On peut donc voir $\mathcal{F}_n$ comme un foncteur de
$\Rep^a(G)$ vers la catégorie des cristaux sur $(X_0/S_n)_{\text{cris}}$.
Avec cette interprétation les $\mathcal{H}^i(\mathcal{A}^s_n)$
s'identifient aux
$R^if_{s,\text{cris}*}\mathcal{O}_{\mathcal{A}^s_0/S_n}$.

\begin{lemm}
Pour tout $V \in \text{Rep}^{a}(G)$, les faisceaux $\mathcal{F}(V)$  et
$\mathcal{F}_n(V)$ sont localement libres sur $X$ et $X_n$
respectivement.
\end{lemm}

En effet c'est le cas pour les $\mathcal{H}^n(\mathcal{A}^s)$.

\subsection{Prolongement des cristaux}

L'objectif est de construire un foncteur $\overline{\mathcal{F}}$ de
$\Rep^a(G)$ vers l'ensemble des fibrés localement libres munis d'une
connexion à pôles logarithmiques le long de $\overline{X}\setminus X$
intégrable et quasi-nilpotente, qui prolonge $\mathcal{F}$.

\subsubsection{Unicité du prolongement}
\label{unicite-prolongement-connexion}

\begin{lemm}
Soit $\mathcal{E}$ un fibré localement libre sur $X$ muni d'une connexion
intégrable et quasi-nilpotente. 

S'il existe un prolongement de $\mathcal{E}$ en un fibré localement libre
sur $\overline{X}$ muni d'une connexion à pôles logarithmiques le long de
$\overline{X}\setminus X$ intégrable et quasi-nilpotente, alors il est
unique, et de plus tout prolongement de $\mathcal{E}$ en un fibré
cohérent sur $\overline{X}$ muni d'une connexion ayant les mêmes
propriétés est aussi localement libre (et donc égal au prolongement
précédent).

De plus, si $\mathcal{E}_1$ et $\mathcal{E}_2$ sont deux tels faisceaux
admettant des prolongements, et $u : \mathcal{E}_1 \to \mathcal{E}_2$ est
un morphisme horizontal, $u$ admet un unique prolongement horizontal
entre les prolongements des faisceaux.
\end{lemm}

\begin{proof}
Soit $\overline{\mathcal{E}}$ un tel prolongement cohérent. Regardons
tout d'abord $\mathcal{E}_K$. D'après \cite{D70}, il existe au plus un
prolongement de $\mathcal{E}_K$ en un fibré muni d'une connexion à pôles
logarithmiques, qui est l'extension canonique de Deligne, ce prolongement
est donc nécessairement $\overline{\mathcal{E}}_K$.

Le faisceau $\overline{\mathcal{E}}$ est alors uniquement déterminé.  En
effet, notons $X'$ la réunion de $X$ et $\overline{X}_K$ dans
$\overline{X}$, $j$ l'inclusion de $X'$ dans $\overline{X}$, et
$\mathcal{E}'$ le faisceau qui coïncide avec $\mathcal{E}$ sur $X$ et
avec $\overline{\mathcal{E}}_K$ sur $\overline{X}_K$.  Alors pour des
raisons de codimension, et le faisceau $\overline{\mathcal{E}}$ étant
cohérent, $\overline{\mathcal{E}} = j_*(\mathcal{E}')$.  En particulier,
tous les prolongements cohérents munis de connexion sont égaux en tant
que faisceaux, donc si l'un est localement libre, tous le sont.

Enfin, $\overline{\mathcal{E}}$ étant localement libre, sa connexion est
entièrement déterminée par sa restriction à $\overline{\mathcal{E}}_K$.

Pour l'existence du prolongement des morphismes, cela provient de la
fonctorialité de l'extension canonique de Deligne.
\end{proof}

\begin{coro}
Le fibré $\mathcal{H}^i(\overline{\mathcal{A}^s})$, muni de la connexion
de Gauss-Manin, ne dépend pas du choix de la compactification
$\overline{\mathcal{A}^s}$. De plus, étant données deux compactifications
$\overline{\mathcal{A}}$ de $\mathcal{A}$ et $\overline{\mathcal{A}^s}$
de $\mathcal{A}^s$, pour tout $i$,
$\mathcal{H}^i(\overline{\mathcal{A}^s})$ et
$\wedge^i\mathcal{H}^1(\overline{\mathcal{A}})^s$ sont égaux comme
sous-faisceaux de $(X\to \overline{X})_*\mathcal{H}^i(\mathcal{A}^s)$
munis d'une connexion à pôles logarithmiques.  
\end{coro}

\begin{proof}
En effet, il suffit pour pouvoir appliquer le lemme précédent de vérifier
que les $\wedge^i\mathcal{H}^1(\overline{\mathcal{A}})^s$ sont localement
libres, il suffit donc de voir que
$\mathcal{H}^1(\overline{\mathcal{A}})$ est localement libre.  Cela se
déduit des résultats de \cite{Ill90}, qu'on peut appliquer car on a
supposé que $\dim_X\mathcal{A} < p$. Le cas général se déduit de
l'identité précédente.
\end{proof}

On notera $\overline{\mathcal{H}}^i(\mathcal{A}^s)$ ce faisceau à connexion.

\begin{lemm}
$\mathcal{H}^i(\overline{\mathcal{A}^s}_n)$ est localement libre sur
$\overline{X}_n$, et ne dépend pas de la compactification
$\overline{\mathcal{A}^s}$.  
\end{lemm}

\begin{proof}
En effet, le faisceau $\overline{\mathcal{H}}^i(\mathcal{A}^s_n)$ est
obtenu à partir de $\overline{\mathcal{H}}^i(\mathcal{A}^s)$ par
changement de base.  
\end{proof} 
On notera dans la suite
$\overline{\mathcal{H}}^i(\mathcal{A}^s_n)$ pour ce faisceau.

\subsubsection{Prolongement de l'action de $E_s$}
\label{existence-prolongement-morphismes-faisceaux}

Il s'agit maintenant de prolonger le morphisme $E_s
\stackrel{a_{\text{cris}}}{\to} \End (\mathcal{H}^{\point}(\mathcal{A}^s))$ en un
morphisme de $\mathbb{Z}_{(p)}$-algèbres 
$E_s \stackrel{a_{\text{log-cris}}}{\to} \End
(\overline{\mathcal{H}}^{\point}(\mathcal{A}^s))$.

La restriction $\text{res} : \End
(\overline{\mathcal{H}}^{\point}(\mathcal{A}^s)) \to \End
(\mathcal{H}^{\point}(\mathcal{A}^s))$ est injective. Pour construire
$a_{\text{log-cris}}$, il suffit donc de vérifier que l'image de
$a_{\text{cris}}$ est contenue dans l'image de $\text{res}$. Comme
$a_{\text{cris}}$ est un morphisme de $\mathbb{Z}_{(p)}$-algèbres, il
suffit de vérifier que l'image par $a_{\text{cris}}$ d'une partie
génératrice de $E_s$ est contenue dans l'image de $\text{res}$. Il s'agit
donc de vérifier que l'action des éléments de $\mathcal{G}$ sur les
$\mathcal{H}^{\point}(\mathcal{A}^s)$ se prolonge en une action sur les
$\overline{\mathcal{H}}^{\point}(\mathcal{A}^s)$.

Soit $u \in \mathcal{G}$. Supposons d'abord que $u$ soit une matrice, ou
(dans le cas unitaire) un élément de $\mathcal{O}_B$. Alors $u$ agit sur
$\mathcal{A}^s$ par une isogénie. D'après la propriété \ref{hyp2}, il
existe donc deux compactifications $\overline{\mathcal{A}^s}_1$ et
$\overline{\mathcal{A}^s}_2$, et un morphisme $u' :
\overline{\mathcal{A}^s}_1 \to \overline{\mathcal{A}^s}_2$ prolongeant
l'action de $u$. Alors $u'$ fournit l'élément de $\End
(\overline{\mathcal{H}}^{\point}(\mathcal{A}^s))$ voulu.  Supposons
maintenant qu'on est dans le cas Siegel et que $u$ est de la forme
$\theta_{i,j}$.  Il suffit de montrer que les morphismes
$a_{\text{cris}}(\varphi_{i,j})$ et $a_{\text{cris}}(\psi_{i,j})$ se
prolongent en éléments de $\End
(\overline{\mathcal{H}}^{\point}(\mathcal{A}^s))$.  Il existe, d'après la
propriété \ref{hyp3}, une compactification $\overline{\mathcal{A}^s}$
telle que le faisceau $(\mathcal{O}(2)^{\otimes a}\otimes
\mathcal{P}_{i,j})^{\otimes b}$ se prolonge en un faisceau $\mathcal{L}$
sur $\overline{\mathcal{A}^s}$, avec $a$ et $b$ premiers à $p$.  Il
existe aussi une compactification $\overline{\mathcal{A}^s}'$ telle que
le faisceau $\mathcal{O}(2)^{\otimes c}$ se prolonge en un faisceau
$\mathcal{L}'$ sur $\overline{\mathcal{A}^s}'$, avec $c$ premier à $p$.
Alors l'action de $\frac{1}{b}(c_\mathcal{L}))$ est dans $\End
(\overline{\mathcal{H}}^{\point}(\mathcal{A}^s))$, l'action de
$-\frac{a}{bc}(c_1(\mathcal{L}'))$ aussi, et l'action de
$\frac{1}{b}(c_1(\mathcal{L})) - \frac{a}{bc}(c_1(\mathcal{L}'))$
prolonge celle de $a_{\text{cris}}(\varphi_{i,j})$.  Pour
$a_{\text{cris}}(\psi_{i,j})$, on fait le même raisonnement, en utilisant
la dualité de Poincaré.

On note encore $a_{\text{log-cris}}$ pour le morphisme $E_s \to
\End(\overline{\mathcal{H}}^{\point}(\mathcal{A}^s_n))$
obtenu par réduction.

\subsubsection{Définition de $\overline{\mathcal{F}}$}

Si $V \in \Rep^a(G)$, on veut définir $\overline{\mathcal{F}}(V)$ comme le
faisceau localement libre muni d'une connexion à pôles logarithmiques
intégrable et quasi-nilpotente sur $\overline{X}$ prolongeant
$\mathcal{F}(V)$.

Au vu du paragraphe \ref{unicite-prolongement-connexion}, il suffit de
montrer l'existence de ce prolongement, son unicité et le fait que la
construction est fonctorielle étant alors automatiques.  Soit $V =
q(\wedge^{\point}\mathcal{V}_0^s)$, o\`u $q$ est un projecteur de $E_s$.
Il suffit de poser $\overline{\mathcal{F}}(V) =
a_{\text{log-cris}}(q)(\overline{\mathcal{H}}^{\point}(\mathcal{A}^s))$.

On note $\overline{\mathcal{F}}_n(V)$ la réduction modulo $\varpi^{n+1}$
de $\overline{\mathcal{F}}(V)$. Munissons $S_n$ de la log-structure
triviale, et $\overline{X}_n$ de la log-structure provenant du diviseur à
croisements normaux $(\overline{X}\setminus X)_n$. Alors
$\overline{\mathcal{F}}_n$ définit un foncteur de $\Rep^a(G)$ vers la
catégorie des cristaux sur
$(\overline{X}_0/S_n)_{\text{cris}}^{\text{log}}$, en effet cette
catégorie est équivalente à celle des
$\mathcal{O}_{\overline{X}_n}$-modules munis d'une connexion à pôles
logarithmiques intégrable et quasi-nilpotente, $\overline{X}_n$ étant un
relèvement de $\overline{X}_0$ log-lisse sur $S_n$.

\section{Structures sur les groupes de cohomologie}

\subsection{Cas étale}

Notons $\overline{K}$ la clôture algébrique de $K$, et $\Gamma =
\text{Gal}(\overline{K}/K)$. Le groupe
$H^m_{\text{ét}}(X_{\overline{K}},\mathbb{F}_n(V))$ est naturellement
muni d'une action de $\Gamma$ car le faisceau $\mathbb{F}_n(V)$ est
défini sur $X_K$. 

On aura besoin du lemme suivant pour comparer l'action de Galois sur
$H^m_{\text{ét}}(X_{\overline{K}},\mathbb{F}_n(V))$
et sur la cohomologie de $\mathcal{A}^s_{\overline{K}}$ : 

\begin{lemm}
\label{decoupage-cohomologie}
Soit $f : Z \to T$ un morphisme de schémas, $F$ un faisceau constant sur
$Z$ ($\mathbb{Z}/p^n\mathbb{Z}$ ou $\mathbb{Z}_p$).
Soit $q$ agissant sur $\mathcal{H}^{\point}(Z) =
R^{\point}f_*F$ et sur $H^{\point}(Z,F)$ de fa\,con compatible avec la
suite spectrale de Leray.
On suppose que $q$ agit comme un projecteur, dont
l'image est entièrement contenue dans $\mathcal{H}^s(Z)$. Notons $V =
q\mathcal{H}^{\point}(Z)$, alors pour tout $m$ on a $H^{m}(T,V) =
qH^{m+s}(Z,F)$.
\end{lemm}

\begin{proof}
En effet considérons la suite spectrale de Leray pour calculer la
cohomologie de $F$ sur $Z$. On lui applique $q$, on obtient toujours une
suite spectrale convergente car $q$ est un projecteur. D'autre part
$qH^m(T,\mathcal{H}^i(Z)) = H^m(T,q\mathcal{H}^i(Z))$, toujours parce que
$q$ est un projecteur. La suite spectrale obtenue a une seule colonne non
nulle, dont les termes sont les $H^m(T,q\mathcal{H}^s(Z))$, et aboutit à
$qH^{m+s}(Z,F)$, d'o\`u le résultat.
\end{proof}

Notons encore $a_{\text{ét}}$ les morphismes $E_s \to
\End(H^{\point}_{\text{ét}}(\mathcal{A}^s_{\overline{K}},\mathbb{Z}_p))$
et $E_s \to
\End(H^{\point}_{\text{ét}}(\mathcal{A}^s_{\overline{K}},\mathbb{Z}/p^n\mathbb{Z}))$.
Supposons que $V = q(\wedge^{\point}\mathcal{V}_0^s)$, $q$ étant un
projecteur homogène de degré $t$, alors on a
$H^m_{\text{ét}}(X_{\overline{K}},\mathbb{F}_n(V)) =
a_{\text{ét}}(q)H^{m+t}_{\text{ét}}(\mathcal{A}^s_{\overline{K}},\mathbb{Z}/p^n\mathbb{Z})$.
Comme $E_s$ agit par des correspondances algébriques définies sur $K$ sur
la cohomologie de $\mathcal{A}^s$, l'action de $\Gamma$ commute à
l'action de $E_s$, et la structure galoisienne obtenue sur
$H^m_{\text{ét}}(X_{\overline{K}},\mathbb{F}_n(V))$ est compatible à
celle sur la cohomologie de $\mathcal{A}^s_{\overline{K}}$.

\subsection{Cas cristallin}
\subsubsection{Les modules de Fontaine-Laffaille}

On note $\underline{MF}^{f}_{\text{tor}}$ la catégorie suivante. Les
objets sont les ${\mathcal{O}}$-modules $M$ de longueur finie, muni d'une
filtration $\text{Fil}^iM$ décroissante, telle que $\text{Fil}^0M = M$ et
$\text{Fil}^{p-1}M = 0$, et pour tout $i$, un $\phi_i : \text{Fil}^iM \to
M$ ${\mathcal{O}}$-semi-linéaire, vérifiant
${\phi_i}_{|\text{Fil}^{i+1}M} = p\phi_{i+1}$, et $\sum_i\im \phi_i = M$.
Les morphismes respectent la filtration et commutent aux $\phi_i$.

\subsubsection{La catégorie $MF(\phi)$}

On introduit la catégorie $MF(\phi)$ des $K$-espaces vectoriels munis
d'une filtration décroissante et d'un Frobenius. Les objets sont les
$K$-espaces vectoriels de dimension finie $M$, muni d'une filtration
décroissante $\text{Fil}^iM$, et de l'action d'un Frobenius $\phi$
semi-linéaire par rapport au Frobenius $\sigma$ de $K$. Les morphismes doivent
commuter au Frobenius, et respecter la filtration.

\subsubsection{Calculs dans un cas particulier}
\label{cas-particulier}

On se place dans le cas suivant : on a un log-schéma $Z$ qui est propre,
et lisse sur $S$ muni de la log-structure triviale. 

Soit $n$ un entier positif.
On note $S_n = \spec {\mathcal{O}}_n$, muni de la log-structure triviale.

Si $(Z,M)$ est un schéma sur $\spec {\mathcal{O}}$, on note $(Z_n,M_n)$
le changement de base à $S_n$.  Si $\mathcal{E}$ est un cristal sur le
site $((Z_m,M_m)/S_n))_{\text{cris}}^{\text{log}}$, on notera
$H^i_{\text{cris}}((Z_m,M_m)/S_n,\mathcal{E})$ pour
$H^i(((Z_m,M_m)/S_n))_{\text{cris}}^{\text{log}},\mathcal{E})$, et
$H^i_{\text{cris}}((Z_m,M_m)/S_n)$ pour
$H^i_{\text{cris}}((Z_m,M_m)/S_n,\mathcal{O}_{Z_m/S_n})$. On omettra la
mention de la log-structure si cela ne cause pas de confusion.
Remarquons que pour tout $m \leq n$, les
$H^i_{\text{cris}}((Z_m,M_m)/S_n)$ ne dépendent pas de $m$, on notera
$H^i_{\text{cris}}((Z,M)/S_n)$ leur valeur commune. Enfin on note
$H^i_{\text{cris}}((Z,M)/S) = \lpro{n}{ H^i_{\text{cris}}((Z,M)/S_n)}$.

On a les deux résultats suivants :

\begin{prop}
Pour tout $0\leq i\leq p-2$, $H^i_{\text{cris}}((Z,M)/S_n)$ est un module de
Fontaine-Laffaille. Pour tout $i\geq 0$,
$H^i_{\text{cris}}((Z,M)/S)\otimes K$ est un élément de $MF(\phi)$.
\end{prop}

\begin{proof}
La preuve de la première partie de la proposition est 
identique à celle de l'article \cite{FonMes87}, qui traite le cas o\`u $Z$
est muni de la log-structure triviale.

Notons $S_n'$ le log-schéma dont le schéma sous-jacent est le même que
$S_n$, et dont la log-structure provient de $\mathbb{N} \to
\mathcal{O}_n, 1 \mapsto 0$.  Il s'agit de la même log-structure que
celle définie dans \cite{HyoKat94}, paragraphe 3.4. Notons $(Z',M')$ le
log-schéma déduit de $(Z,M)$ par le changement de base $S_n' \to S_n$.
Alors $H^i_{\text{cris}}((Z,M)/S_n)$ et $H^i_{\text{cris}}((Z',M')/S_n')$
sont canoniquement isomorphes pour tout $i$. Les
$H^i_{\text{cris}}((Z',M')/S_n')$ sont munis d'un Frobenius et d'un
opérateur de monodromie, définis dans \cite{HyoKat94}, paragraphe 3.
L'opérateur de monodromie est ici nul, $(Z',M')$ provenant par changement
de base de $(Z,M)$ qui est log-lisse sur $S_n$.
$H^i_{\text{cris}}((Z',M')/S')\otimes K$, et donc aussi
$H^i_{\text{cris}}((Z,M)/S)\otimes K$ est ainsi naturellement muni d'une
structure d'élément de $MF(\phi)$.  \end{proof}

\subsubsection{Action des endomorphismes sur le prolongement des
cristaux}

Notons  $\overline{H}^i(\mathcal{A}^s/S_n)$ la limite des
$H^i_{\text{cris}}(\overline{\mathcal{A}^s}_n/S_n)$, pour les 
$\overline{\mathcal{A}^s}$ de notre famille de compactifications. Alors : 

\begin{prop}
Le morphisme $H^i_{\text{cris}}(\overline{\mathcal{A}^s}_n/S_n)\to\overline{H}^i(\mathcal{A}^s/S_n) $ 
est un isomorphisme pour toute compactification $\overline{\mathcal{A}^s}$ et
pour tout $n$.
\end{prop}

\begin{proof}
Il suffit pour cela de voir que tout morphisme entre compactifications
qui est l'identité sur $\mathcal{A}^s$ induit un isomorphisme entre les
groupes de cohomologie.

Considérons la suite spectrale de Leray :
$E^{i,j}_2 = 
H^i(\overline{X}_n/S_n,R^jf_{s\text{cris}*}\mathcal{O}_{\overline{\mathcal{A}^s}_n/S_n}) 
\Rightarrow
H^{i+j}(\overline{\mathcal{A}^s}_n/S_n)$

Un morphisme entre deux compactifications induit un morphisme de suites
spectrales, qui est un isomorphisme sur les $E^{i,j}_2$, donc aussi sur
l'aboutissement. 
\end{proof}

On cherche à définir un morphisme de $\mathbb{Z}_{(p)}$-algèbres $E_s \to
\End(\overline{H}^i(\mathcal{A}^s/S_n))$.

On a la suite spectrale suivante, qu'on appellera encore suite spectrale
de Leray, qui provient de n'importe quelle compactification
$\overline{\mathcal{A}^s}$ de $\mathcal{A}^s$ : 
$$ 
E^{i,j}_2 =
H^i(\overline{X}_n/S_n,\overline{\mathcal{H}}^j(\mathcal{A}^s_n))
\Rightarrow \overline{H}^{i+j}(\mathcal{A}^s/S_n) 
$$

\begin{prop}
Il existe un unique morphisme $a_{\text{log-cris}}$ de
$\mathbb{Z}_{(p)}$-algèbres $E_s \to
\End(\overline{H}^i(\mathcal{A}^s/S_n))$ tel que l'action de $E_s$ sur
les $\overline{\mathcal{H}}^j(\mathcal{A}^s_n)$ et sur les
$\overline{H}^i(\mathcal{A}^s/S_n)$ donnée par $a_{\text{log-cris}}$ soit
compatible à la suite spectrale de Leray.  
\end{prop}

\begin{proof}
On commence par définir l'image de l'ensemble $\mathcal{G}$ des éléments
géométriques de $E_s$, et on montre ensuite que l'on peut prolonger en un
morphisme de $\mathbb{Z}_{(p)}$-algèbres.

Pour définir l'image d'un élément de $\mathcal{G}$, on fait exactement
comme dans le paragraphe
\ref{existence-prolongement-morphismes-faisceaux}.
Il faut voir que le choix fait est unique. 
Cela provient du fait que l'action de $u \in \mathcal{G}$ sur les termes
$E^{i,j}_2$ de la suite spectrale ne dépend pas des choix faits, comme
expliqué en \ref{existence-prolongement-morphismes-faisceaux}, et de la
compatibilité de l'action du prolongement à la suite spectrale de Leray. 

Il reste à voir que l'action des éléments de $\mathcal{G}$ se prolonge en
un morphisme d'algèbres $E_s \to
\End(\overline{H}^i(\mathcal{A}^s/S_n))$. Cela provient encore une fois
de la compatibilité avec la suite spectrale de Leray, et du fait que $E_s
\to \End(\overline{\mathcal{H}}^{\point}(\mathcal{A}^s/S_n))$ est un
morphisme d'algèbres.  
\end{proof}

\subsubsection{La cohomologie des cristaux}
\label{coho-cristaux}

Posons $H^i_{\text{cris}}(\overline{X}/S,\overline{\mathcal{F}}(V))=
\lpro{n}H^i_{\text{cris}}(\overline{X}_n/S_n,\overline{\mathcal{F}}_n(V))$.
On a le résultat suivant :

\begin{prop}
Pour tout $V\in \Rep^a(G)$, pour tout $i$,
$H^i_{\text{cris}}(\overline{X}/S,\overline{\mathcal{F}}(V))\otimes K$ est un élément de $MF(\phi)$.

Pour tout $V\in \Rep^a(G)$ homogène de degré $t$, pour tout $i$ tel que
$i+t\leq p-2$, pour tout $n$,
$H^i_{\text{cris}}(\overline{X}_n/S_n,\overline{\mathcal{F}}_n(V))$ est un
élément de $\underline{MF}^{f}_{\text{tor}}$. 
\end{prop}

\begin{proof}
En effet, soit $V = q(\wedge^{\point}\mathcal{V}_0^s)$, avec $q$
homogène de degré $t$, on a alors de fa\,con similaire au lemme
\ref{decoupage-cohomologie} les égalités 
$H^i_{\text{cris}}(\overline{X}/S,\overline{\mathcal{F}}(V))\otimes K =
a_{\text{cris}}(q)\overline{H}^{i+t}_{\text{log-cris}}(\mathcal{A}^s/S)\otimes K$ 
et 
$H^i_{\text{cris}}(\overline{X}_n/S_n,\overline{\mathcal{F}}_n(V)) =
a_{\text{cris}}(q)\overline{H}^{i+t}_{\text{cris}}(\mathcal{A}^s_n/S_n)$.
Il reste à voir que l'action de $E_s$ par $a_{\text{log-cris}}$ respecte
les structures d'élément de $MF(\phi)$, et que l'action de $E(A)_s$
respecte les structures de module de Fontaine-Laffaille, ce qui se voir
sur les éléments géométriques.
\end{proof}

\section{Théorème de comparaison}
\label{comparaison}

\subsection{Le cas des faisceaux constants}
\label{comparaison-constant}

Notons $\Rep_{\mathbb{Z}_p}(\Gamma)$ la catégorie des
$\mathbb{Z}_p$-représentations de type fini de $\Gamma =
\text{Gal}(\overline{K}/K)$. Nous avons un foncteur contravariant et
pleinement fidèle : $V_{\text{cris}} : \underline{MF}^{f}_{\text{tor}}
\to \Rep_{\mathbb{Z}_p}(\Gamma)$ qui est défini par $V_{\text{cris}}(M) =
\Hom(M,A_{\text{cris},\infty})$.  L'anneau $A_{\text{cris}}$ est défini
dans \cite{Bre96}, 6.3, et $A_{\text{cris},\infty} = A_{\text{cris}}
\otimes \mathbb{Q}_p/\mathbb{Z}_p$. L'anneau $A_{\text{cris}}$ est muni
d'une filtration décroissante et d'un Frobenius, et les homomorphismes
que l'on considère doivent être compatibles à la filtration et à l'action
du Frobenius.

Soit $Z$ un schéma propre et lisse sur $\spec \mathcal{O}_K$, et $D$ un
diviseur à croisements normaux relatifs de $Z$, $U$ l'ouvert
complémentaire de $D$.  On munit $Z$ de la log-structure $M$ définie par
le diviseur $D$, et $\spec \mathcal{O}_K$ de la log-structure triviale.

\begin{prop}
\label{constant-torsion}
Pour $0 \leq m \leq p-2$, on a un isomorphisme canonique
compatible à l'action de Galois :
$$
V_{\text{cris}}(H^m_{\text{cris}}((Z,M)/S_n)) =
\chk{H^m_{\text{ét}}(U_{\overline{K}},\mathbb{Z}/p^n\mathbb{Z})}
$$
\end{prop}

\begin{prop}
\label{constant-rationnel}
Pour tout $m$, il existe un isomorphisme canonique qui respecte l'action
de $\Gamma$, la filtration et le Frobenius :
$$
\gamma_m :
B_{\text{cris}}\otimes_{\mathcal{O}}H^m_{\text{cris}}((Z,M)/S)
\isom 
B_{\text{cris}}\otimes_{\mathbb{Q}_p}H^m_{\text{ét}}(U_{\overline{K}},\mathbb{Q}_p)
$$
\end{prop}

\begin{proof}
Le résultat de la proposition \ref{constant-torsion} provient de travaux
de Breuil (\cite{Bre98}) et Tsuji (\cite{Tsu00}). Ces résultats
s'appliquent dans un cadre beaucoup plus général que celui considéré ici,
et décrivent une comparaison entre la cohomologie étale de
$U_{\overline{K}}$ et la cohomologie de $(Z',M')/E_n$. Ici $(Z',M')$ est
obtenu comme dans le paragraphe \ref{cas-particulier} par changement de
base de $(Z,M)$ de $S_n$ à $S_n'$.  $E_n$ est le log-schéma dont le
schéma sous-jacent est $\spec \mathcal{O}_n\langle u\rangle$, l'enveloppe
à puissances divisées de l'algèbre $\mathcal{O}_n[u]$ des polynômes en
l'indéterminée $u$, muni de la log-structure associée à $\mathbb{N} \to
\mathcal{O}_n\langle u\rangle$, $1 \mapsto u$. Dans notre cas
particulier, on a une relation simple entre la cohomologie de
$(Z',M')/E_n$ et celle de $(Z',M')/S_n'$, donnée par
$H^m_{\text{cris}}((Z',M')/E_n) = \mathcal{O}_n\langle u\rangle \otimes
H^m_{\text{cris}}((Z',M')/S_n')$, qui nous permet d'obtenir le résultat
de la proposition \ref{constant-torsion}.

Pour la version rationnelle \ref{constant-rationnel}, le résultat provient de
résultats de Tsuji (\cite{Tsu99}, voir aussi \cite{YamI}). 
Comme dans le cas de torsion, la situation se simplifie par rapport au
cas général, du fait qu'ici la monodromie agissant sur
$H^m_{\text{cris}}((Z,M)/S)\otimes K$ est nulle.
\end{proof}

\subsection{Les théorèmes}

\begin{theo}
\label{principal_torsion}
Dans le cas unitaire, soit $V \in \Rep^a(G)$, et $m$ tel que $m+t(V) \leq
p-2$.  $H^m_{\text{ét}}(X_{\overline{K}},\mathbb{F}_n(V))$ est muni d'une
action du groupe de Galois $\Gamma$,
${H}^m_{\text{cris}}(\overline{X}_n/S_n,\overline{\mathcal{F}}_n(V))$ est
muni d'une structure de module de Fontaine-Laffaille, et on a un
isomorphisme : 
$$
V_{\text{cris}}({H}^m_{\text{cris}}(\overline{X}_n/S_n,\overline{\mathcal{F}}_n(V))
)= \chk{H^m_{\text{ét}}(X_{\overline{K}},\mathbb{F}_n(V))} 
$$ 
\end{theo}

Notons qu'on peut déduire de ce théorème comme dans l'article
\cite{Bre98}, paragraphe 4.2, une comparaison entre les parties de
torsion de $\lpro{n} H^m_{\text{ét}}(X_{\overline{K}},\mathbb{F}_n(V))$
et de $\lpro{n}
{H}^m_{\text{cris}}(\overline{X}_n/S_n,\overline{\mathcal{F}}_n(V))$,
ainsi qu'une comparaison entre leurs parties libres.

\begin{theo}
\label{principal_rationnel}
Dans le cas unitaire et Siegel, soit $V \in
\Rep^a(G)$, il existe un isomorphisme 
$$\gamma : 
B_{\text{cris}}\otimes_{\mathcal{O}}{H}^m_{\text{log-cris}}(\overline{X}/S,\overline{\mathcal{F}}(V))
\isom 
B_{\text{cris}}\otimes_{\mathbb{Z}_p}H^m_{\text{ét}}(X_{\overline{K}},\mathbb{F}(V))
$$
\end{theo}

Le point essentiel de la preuve dans les deux cas est le résultat suivant :

\begin{lemm}
\label{compatibilite-actions}
Soit $u \in E(A)_s$. $u$ agit sur $\overline{H}^m(\mathcal{A}^s/S_n)$ et
sur
$H^m_{\text{ét}}(\mathcal{A}^s_{\overline{K}},\mathbb{Z}/p^n\mathbb{Z})$
($m \leq p-2$) de fa\,con compatible avec l'isomorphisme
$V_{\text{cris}}$. Soit $u \in E_s$, $u$ agit sur
$\overline{H}^m(\mathcal{A}^s/S)\otimes\mathbb{Q}$ et sur
$H^m_{\text{ét}}(\mathcal{A}^s_{\overline{K}},\mathbb{Q}_p)$ de fa\,con
compatible avec l'isomorphisme $\gamma_m$ du théorème
\ref{constant-rationnel}.  
\end{lemm}

\begin{proof}
Il suffit de montrer la compatibilité des actions pour l'ensemble des
éléments géométriques $\mathcal{G}$ de $E_s$, puisqu'ils engendrent $E_s$
comme $\mathbb{Z}_{(p)}$-algèbre.

Soit $u$ un élément de $E_s$ provenant d'une matrice de déterminant non
nul, ou d'un élément non nul de $\mathcal{O}_B$. Son action sur la
cohomologie provient d'une isogénie de $\mathcal{A}^s$, qu'on notera
encore $u$.  D'après la propriété \ref{hyp2} $u$ se prolonge en un
morphisme entre deux compactifications $u : \overline{\mathcal{A}^s}_1
\to \overline{\mathcal{A}^s}_2$.  D'o\`u par fonctorialité de
$V_{\text{cris}}$, $V_{\text{cris}}(u :
H^m_{\text{cris}}((\overline{\mathcal{A}^s}_2)_n/S_n)\to
H^m_{\text{cris}}((\overline{\mathcal{A}^s}_1)_n/S_n) = (\chk{u} :
\chk{H^m_{\text{ét}}(\mathcal{A}^s_{\overline{K}},\mathbb{Z}/p^n\mathbb{Z})}
\to
\chk{H^m_{\text{ét}}(\mathcal{A}^s_{\overline{K}},\mathbb{Z}/p^n\mathbb{Z})})$,
ce qui est bien la compatibilité voulue. De m\^eme, on a aussi la
compatibilité pour l'action sur
$H^m_{\text{ét}}(U_{\overline{K}},\mathbb{Q}_p)$ et
$H^m_{\text{cris}}((X,M)/S)\otimes K$.

Dans le cas Siegel, il faut aussi considérer les éléments de la forme
$\theta_{i,j}$. Il s'agit donc de voir que l'action des $\varphi_{i,j}$
et des $\psi_{i,j}$ sur $H^m_{\text{ét}}(U_{\overline{K}},\mathbb{Q}_p)$
et \mbox{$H^m_{\text{cris}}((X,M)/S)\otimes K$} est compatible. Cela
provient du fait que l'isomorphisme de comparaison
\ref{constant-rationnel} fait correspondre les classes de Chern
(\cite{Tsu99}) et est compatible aux structures produit sur les groupes
de cohomologie et à la dualité de Poincaré (\cite{YamI}).  \end{proof}

\begin{proof}[Démonstration des théorèmes \ref{principal_torsion} et
\ref{principal_rationnel}]
Montrons le théorème \ref{principal_torsion}.  Soit $V \in \Rep^a(G)$, on
peut supposer qu'il existe un entier $s$, et un projecteur $q$ dans
$E(A)_s$ de degré $t$, tels que $V = q(\wedge^{\point}\mathcal{V}_0^s)$.
Le théorème de comparaison s'applique car \mbox{$m+t \leq p-2$}, et nous
donne un isomorphisme
$V_{\text{cris}}(\overline{H}^{m+t}(\mathcal{A}^s/S_n)) =
\chk{H^{m+t}_{\text{ét}}(\mathcal{A}^s_{\overline{K}},\mathbb{Z}/p^n\mathbb{Z})}$.

Appliquons $q$ : comme l'action de $E(A)_s$ commute à $V_{\text{cris}}$
on a donc un isomorphisme :
$V_{\text{cris}}(a_{\text{log-cris}}(q)\overline{H}^{m+t}(\mathcal{A}^s/S_n))
=
\chk{(a_{\text{ét}}(q)H^{m+t}_{\text{ét}}(\mathcal{A}^s_{\overline{K}},\mathbb{Z}/p^n\mathbb{Z}))}$.

D'après le lemme \ref{decoupage-cohomologie}, cela donne :
$V_{\text{cris}}({H}^m_{\text{log-cris}}(\overline{X},\overline{\mathcal{F}}_n(V)))
= \chk{H^m_{\text{ét}}(X_{\overline{K}},\mathbb{F}_n(V))}$.

La preuve du théorème \ref{principal_rationnel} est identique.
\end{proof}

\begin{rema}
On voit appara\^itre dans le lemme \ref{compatibilite-actions} le point
qui explique pourquoi on n'a pas de résultats de comparaison prenant en
compte la torsion pour le cas Siegel : la compatibilité de l'isomorphisme
de comparaison à coefficients constants avec la dualité de Poincaré et
les structures produits n'est actuellement montrée que dans le cas
rationnel (même s'il est vraisemblable qu'elle soit vraie aussi dans le
cas de torsion, en introduisant des limitations sur le degré des groupes
de cohomologie considérés).  Enfin on peut remarquer que si on se limite
aux représentations qui peuvent être obtenues à l'aide uniquement des
éléments de $E_s$ provenant de $M_s(\mathbb{Z})$, on peut prendre en
compte la torsion pour le cas Siegel.  
\end{rema}

\providecommand{\bysame}{\leavevmode ---\ }
\providecommand{\og}{``}
\providecommand{\fg}{''}
\providecommand{\smfandname}{et}
\providecommand{\smfedsname}{\'eds.}
\providecommand{\smfedname}{\'ed.}
\providecommand{\smfmastersthesisname}{M\'emoire}
\providecommand{\smfphdthesisname}{Th\`ese}

\end{document}